\documentclass[10pt]{amsart}
\usepackage{graphicx} 
\usepackage{amsfonts,amsmath,tikz,mathabx,caption}
\usepackage{amssymb}
\usepackage{pb-diagram}
\usepackage[all]{xy}
\usepackage{hyperref,xcolor}
\usetikzlibrary{decorations.markings}
\usetikzlibrary{shapes}
\usetikzlibrary{backgrounds}
\usepackage{subfig}

  \newtheorem{thm}{Theorem}
\newtheorem{lemma}[thm]{Lemma}

\newtheorem{cor}[thm]{Corollary}
\newtheorem{prop}[thm]{Proposition}

\numberwithin{thm}{section}

\newcommand{\Z}{\mathbb{Z}}
\newcommand{\N}{\mathbb{N}}

\newcommand{\R}{\mathbb{R}}

\newcommand{\Sym}{\textup{Sym}}
\newcommand{\Alt}{\textup{Alt}}

\newcommand{\Supp}{\textup{Supp}}
  
\title[MIF growth]{Finitely generated groups with arbitrarily fast MIF growth}
\author{Yash Lodha}
\address{Department of Mathematics, Purdue University, West Lafayatte, IN, USA}
    \email{ylodha@purdue.edu}
    \urladdr{https://yl7639.wixsite.com/website}
\date{\today}

\begin{document}

\begin{abstract}
A group $G$ is mixed identity free, or MIF, if for each nontrivial $w(x)\in G\star (\mathbb{Z}=\langle x\rangle)$ there is some $g\in G$ such that $w(g)\neq 1$.
Associated to MIF groups is a natural growth function that measures how far in the Cayley graph one witnesses such inequations. 
We adapt a recent construction of the author with Hyde to show that for each function $f:\N\to \N$ there is a finitely generated group $G$ whose MIF growth function grows faster than $f$.
This answers a question first considered by Henry Bradford \cite{BradfordPC}. As a consequence, we also provide the first examples of finitely generated MIF groups that are not selfless, answering a question of Bradford and Sisto in \cite{BradfordSisto2026}.

\end{abstract}

\maketitle

\section{Introduction.}

Mixed identities in groups were first defined in an article of Anashin \cite{Anashin1978} in $1978$. 
A group $G$ is mixed identity free (or MIF), if for each $w(x)\in G\star (\Z=\langle x\rangle)$ such that $w(x)\neq 1$,
there is some $g\in G$ such that $w(g)\neq 1$. 
This property was systematically studied and popularized in an article of Hull and Osin \cite{HullOsin2016}, who established that acylindrically hyperbolic groups with trivial finite radical are MIF.
This has come to be a fundamental property in group theory, and has recently witnessed some surprising applications towards longstanding open problems in $C^*$-algebras (see \cite{AmrutamEtAl2025} and \cite{KunnawalkamElayavalliSchafhauser2025}), as well as striking connections with the Boone-Higman conjecture \cite{BFFHZ25}.

A fundamental observation of Hull and Osin in \cite{HullOsin2016} is the following dichotomy for a highly transitive faithful action of a group $G$ on a set $X$. Either $G$ contains as a normal subgroup $\textup{Alt}_f(X)$, the group of finitely supported alternating permutations of $X$, or $G$ is mixed identity free. 
Here a group is said to be \emph{highly transitive} if it admits a faithful action on $\N$ which is transitive on the set of ordered $n$-tuples of distinct elements of $\N$, for each $n\in \N_{\geq 1}$.
In \cite{HL25}, the author and Hyde provided a general construction of highly transitive and mixed identity free groups which are also amenable.

Given a finite generating set $S$ for for a MIF group $G$, the complexity of $w(x)\in G\star \Z\setminus \{1\}$ is defined as: $$C_G(w(x))=\textup{min}\{||g||_S\mid w(g)\neq 1\}$$
One defines the MIF growth function $\mathcal{M}_G^S:\N_{\geq 1}\to \N$ where: $$\mathcal{M}_G^S(n)=\textup{max}\{C(w(x))\mid ||w(x)||_{S\cup \{x\}}\leq n, w(x)\neq 1\}$$
The growth type is independent of the choice of the generating set.
This growth function was first defined by Henry Bradford in \cite{Bradford2024Lawlessness}.
In this article, we prove the following.

\begin{thm}\label{thm:main}
Let $f:\N\to \N$ be a function. Then there is a finitely generated MIF group $G$ whose MIF growth function eventually dominates $f$ pointwise. Moreover, such a $G$ can be chosen to be $2$-generated and highly transitive.
\end{thm}

Recall that for a finitely generated group $G$ with an ordered fixed finite generating set $S=(s_1,\ldots, s_k)$, the pair $(G,S)$ is called a marked group. 
$\Gamma(G,S)$ is the associated Cayley graph with vertex set $G$ and edge set $\{(g,gs)\mid g\in G,s\in S\}$. This is endowed with the standard path metric, and for $g\in G$, we denote the usual word norm as $|g|_S$. $G$ acts by left multiplication on $\Gamma(G,S)$ by isometries. We denote by $B_n(\Gamma(G,S))$ the closed $n$-ball centered around the identity in the Cayley graph and by $B_n^v(\Gamma(G,S))$ we denote the set of elements in $G$ whose $S$-word norm is at most $n$ (i.e. the vertices in $B_n(\Gamma(G,S))$). Note that here we make a formal distinction between the graph and its underlying vertex set.

A stronger property than MIF, called \emph{selflessness}, was introduced by the authors in \cite{AmrutamEtAl2025} which plays a crucial role in their solution to Leonel Robert’s selflessness problem for reduced group $C^*$-algebras. Let $G$ be a group with a finite generating set $S$. Let $H=G\star \Z$ be endowed with the finite generating set $S\cup \{z\}$ where $\Z=\langle z\rangle$. Then the group $G$ is called selfless if there is a function $f:\N\to \R$ such that $\liminf \sqrt[n]{f(n)}=1$ and such that for each $n\in \N$ there is a homomorphism $\phi_n:G\star Z\to G$ such that the following holds:
\begin{enumerate}
\item $\phi_n\restriction G=id$.
\item $\phi_n$ is injective on $B_n^v(H,S\cup \{z\})$.
\item $\phi_n(B_n^v(H,S\cup \{z\}))\subseteq B_{f(n)}^v(G,S)$.
\end{enumerate}
This does not depend on the choice of finite generating set $S$.
By \cite[Definition~3.1 and Remark~3.2(i)]{BradfordSisto2026},
if $G$ is selfless, then there exists a function $f$ satisfying
$\liminf_{n\to\infty} f(n)^{1/n}=1$ such that
$\mathcal{M}_G^S(n)\leq f(n)$ for every $n$.
Consequently, a finitely generated MIF group whose MIF growth
function admits an eventual exponential lower bound is not selfless.
So our main theorem has the following corollary, which answers Question~3.3 of Bradford and Sisto
\cite{BradfordSisto2026} in the negative.

\begin{cor}
There exist finitely generated MIF groups that are not selfless.
\end{cor}

\section*{Acknowledgement}
The author is supported by the NSF CAREER grant DMS-2552707. The paper is entirely human written, and the proofs are due to the author. AI help was used to construct the bibliography.

\section{The construction}

We view $\Gamma(G,S)$ and $B_n(\Gamma(G,S))$ as marked graphs, i.e. each edge $(g,gs)$ is marked by the label $s$ and oriented with $g$ as the source.
Two marked groups $(G_1,S_1=(s_1,\ldots, s_k)), (G_2,S_2=(t_1,\ldots t_k))$ are said to be \emph{marked isomorphic} if the bijection $S_1\to S_2, s_i\mapsto t_i$ extends to an isomorphism $G_1\to G_2$. We consider the associated equivalence relation among $n$-generated marked groups for a fixed $n$, and define the space: $$\mathcal{G}_n=\{(G,S)\mid |S|=n\}/\sim$$
The space is called the \emph{Grigorchuk space of marked groups} and is endowed with the topology determined by the metric:
$d((G_1,S_1),(G_2,S_2))=2^{-m}$ where $m$ is the largest number such that $B_m(\Gamma(G_1,S_1))\cong B_m(\Gamma(G_2,S_2))$ (marked isomorphic and mapping the identity to the identity).
Here the isomorphism between the $m$-balls is an isomorphism as a marked graph (i.e. an edge labelled by the $i$'th generator of $S_1$ must map to an edge labelled by the $i$'th generator of $S_2$).

All actions of permutation groups in this paper will be right actions.
The following construction is a controlled variation of a special case of the construction by the author with Hyde in \cite{HydeLodha2025}.
However, the proofs below are self contained.
Recall that $\Z$ embeds in $\Sym(\Z)$ via the right regular action. $\Sym(\Z)$ is endowed with the topology of pointwise convergence. 
The same for $\Alt_f(\Z)$, the subgroup consisting of alternating finitely supported permutations.
Given an element $f\in \Sym(\Z)$, the support of the element, or $\Supp(f)$ is the set of points in $\Z$ that are moved by $f$.
We say that $f,g\in \Sym(\Z)$ are \emph{disjoint} if they have disjoint supports. 
Given an infinite collection of pairwise disjoint nontrivial permutations $\{f_{m}\}_{m\in \N}$, we can consider their product $\prod_{m\in \N}f_{m}$ which is the limit $\lim_{n\to \infty}f_1...f_n$ in the topology of pointwise convergence. 
We denote by $\Sym_f(\Z)$ to be the set of finitely supported permutations.
An element $f\in \Sym(\Z)$ is called an \emph{involution} if it is an element of order two, or equivalently, in cycle notation it is a (possibly infinite) product of disjoint transpositions. We denote $[x,y]=x^{-1}y^{-1}xy$ and $y^x=x^{-1}yx$.

The construction is motivated by the following elementary observation.
Consider $\Z=\langle z\rangle$.
We view the natural right action of $\Z$ on itself as an embedding $\Z\to \Sym(\Z)$.
Therefore, $n\cdot z=n+1$ for each $n\in \Z$.
Let $f\in \Sym_f(N)$ be an involution so that $\Supp(f)\subseteq [0,\ldots , m]$ for $m\in \N$.
Consider the marked group $$(H_0=\langle f,z\rangle, K_0=\{f,z\})$$
Let $k\in \N$ and consider the marked graph $B_k(\Gamma(H_0,K_0))$.
Then the following holds:

\begin{lemma}\label{lem:elemob}
There is an $l\in \N$ such that the marked group 
$$(H_1,K_1)\qquad H_1=\langle [z^l,f], z\rangle\qquad  K_1=\{[z^l,f],z\}$$
satisfies that $$B_k(\Gamma(H_0,K_0))\cong B_k(\Gamma(H_1,K_1))\qquad [z^l,f]^2=1$$
\end{lemma}

\begin{proof}
Recall that $f$ is an involution.
For $l$ sufficiently large, the support of $[z^l,f]$ is contained in the union of two blocks $$\{0,\ldots, m\}\cup \{l,\ldots, l+m\}$$
so that $[z^l,f]\restriction \{l,\ldots, l+m\}$ is a conjugate of $[z^l,f]\restriction \{0,\ldots, m\}$ by $z^l$. 
Since the element $z$ acts ``homogeneously" on $\Z$ by translation, $l$ can be chosen to be sufficiently large so that words of the prescribed length in the generating set $\{f, z\}$ satisfy the same relations as those in $\{[z^l,f],z\}$.
\end{proof}

We fix an increasing function $\zeta:\N\to \N$ for the rest of this article. It will be our goal to construct a marked group $(G,S=\{z,\sigma\})$ such that:
$$\mathcal{M}_G^S(n)\geq \zeta(n)\qquad \text{ for sufficiently large } n$$
We also fix the increasing sequences: 
$$e_n=(2^{n+2})!\qquad q_n=\zeta(4e_{n+1})\qquad R_n=2e_n(q_n+1)$$ for the rest of the article, which will play a crucial role in the construction of the groups.

Let $\Z=\langle z\rangle$ be as above and consider the permutation $\sigma_0=(0,1)\in \Sym(\Z)$. 
We define $S_0=\{ z,\sigma_0\}$ and the group $G_0=\langle S_0\rangle$.

We shall construct a sequence $(l_n)_{n\in \N}$ of natural numbers with $l_0=1$.
We define in an inductive fashion $\sigma_n=[z^{l_n},\sigma_{n-1}]$ so that by design, each $\sigma_n$ is an involution and 
$|\Supp(\sigma_n)|=2^{n+1}$.
We define the marked groups $(G_n,S_n)$ where $S_n=(z,\sigma_n)$. 
Assume that we have defined $$(l_m)_{m\in\{0,\ldots ,n\}}\qquad (G_m,S_m=\{z,\sigma_m\})_{m\in \{0,\ldots, n\}}$$
First observe that the following holds.
\begin{lemma}\label{lem:unifMI}
The group $G_n$ satisfies the nontrivial mixed identity $w_n(x)=[x,\sigma_n]^{e_n}$.
\end{lemma}

\begin{proof}
This follows immediately from the fact that $$|\Supp(\sigma_n)|=|\Supp(x^{-1}\sigma_n^{-1}x)|=2^{n+1}$$ for any $x\in G_n$, and hence $|\Supp([x,\sigma_n])|\leq 2^{n+2}$.
\end{proof}
We fix the words 
$$w_n(x)=[x,\sigma_n]^{e_n}\in G_n\star (\Z=\langle x\rangle)$$ for the rest of the article.
Note that these are nontrivial in the respective free product, cyclically reduced, and of length $4e_n$.
Now we must define the number $l_{n+1}$. 
Recall that $e_{n+1}=(2^{n+3})!$ and $q_n=\zeta(4e_{n+1})$. First, we observe the following.

\begin{lemma}\label{lem:estimate}
There is an $l_{n+1}\in \N$ such that $\frac{l_{n+1}}{l_n}\geq 2$ and there is a marked isomorphism: 
$$B_{R_n}(\Gamma(G_{n},S_{n}))\cong  B_{R_n}(\Gamma(G_{n+1},S_{n+1}))$$ where $$S_{n+1}=(z,\sigma_{n+1}=[z^{l_{n+1}},\sigma_{n}])\qquad G_{n+1}=\langle S_{n+1}\rangle$$
\end{lemma}

\begin{proof}
This follows from a direct application of Lemma \ref{lem:elemob}.
\end{proof}




Lemma \ref{lem:estimate} allows us to define the marked group $$(G_{n+1}, S_{n+1}=(z, \sigma_{n+1}=[z^{l_{n+1}},\sigma_{n}]))$$
By design, we have a well defined limit $\sigma=\lim_{n\to \infty} \sigma_n$ in $\Sym(\Z)$, which is also an involution.
We define the marked group $G=\langle \Z,\sigma\rangle$ with the marking $S=\{z,\sigma\}$.
We also define the words $$w_{n,\infty}(x)=[x,\sigma]^{e_n}\in G\star (\Z=\langle x\rangle)$$
and observe that these have $S\cup \{x\}$-length $4e_n$.

\begin{lemma}\label{lem:length}
Let $g\in G$ be such that $|g|_S\leq q_n$.
Then $|w_{n,\infty}(g)|_S\leq R_n$.
\end{lemma}

\begin{proof}
This follows since $w_{n,\infty}(x)=[x,\sigma]^{e_n}$, and upon substituting $x$ by $g$.
\end{proof}

\begin{prop}
For each $n\in \N$ there is a marked isomorphism
$$B_{R_n}(\Gamma(G_{n},S_{n}))\cong  B_{R_n}(\Gamma(G,S))$$
In particular, it follows that $\zeta(4e_{n+1})= q_{n}\leq C_G(w_{n,\infty}(x))$.
\end{prop}

\begin{proof}
The proof of the first part is immediate from the definitions and from the fact that $(R_n)_{n\in \N}$ is an increasing sequence.
We will show that $\zeta(4e_{n+1})= q_{n}\leq c_{G}(w_{n,\infty}(x))$. 
Recall that $$w_n(x)=[x,\sigma_n]^{e_n}\qquad w_{n,\infty}(x)=[x,\sigma]^{e_n}$$
and that $$q_n=\zeta(4e_{n+1})\qquad R_n=2e_n(q_n+1)$$ 
It suffices to show that given any $g\in G$ such that $|g|_S\leq q_n$,
we have that $w_{n,\infty}(g)=1$.
It follows from Lemma \ref{lem:length} that $|w_{n,\infty}(g)|_{S}\leq R_n$.
Let $g_1$ be the $S_n$-word which maps to $g\in G$ under the marked isomorphism $B_{R_n}(\Gamma(G_{n},S_{n}))\cong  B_{R_n}(\Gamma(G,S))$.
Note that $|w_{n}(g_1)|_{S_n}\leq R_n$.
Since $w_n(g_1)=1$ in $G_n$ and since we have a marked isomorphism, we also conclude that $w_{n,\infty}(g)=1$ in $G$.
\end{proof}

This provides the following estimate on the MIF growth.

\begin{prop}
For each $n\in \N,n\geq 4e_1$, we have that $\mathcal{M}_{G}^{S}(n)\geq \zeta(n)$.
\end{prop}

\begin{proof}
First we observe that given a sequence of increasing natural numbers $(c_n)_{n\in \N}$, and two nondecreasing functions $g_1,g_2:\N\to \N$ satisfying that 
$g_1(c_n)\geq g_2(c_{n+1})$, it holds that $g_1(m)\geq g_2(m)$ for each $m\geq c_1$.
This follows since whenever $m\in (c_n,c_{n+1})$, we have $$g_2(m)\leq g_2(c_{n+1})\leq g_1(c_n)\leq g_1(m)$$
We apply this to the sequence $(4e_n)_{n\in \N}$, using the fact that $$\zeta(4e_{n+1})\leq C_G(w_{n,\infty}(x))\leq  \mathcal{M}_{G}^{S}(|w_{n,\infty}(x)|)=\mathcal{M}_{G}^{S}(4e_n)$$
\end{proof}

Our final goal is to show that the group $G$ is highly transitive and mixed identity free. 
The following Lemma serves as a key ingredient.
\begin{lemma}\label{lem:HT}
For each $i,n\in \Z, i<n$, there is an element $\nu\in G$ such that $\nu\restriction (-\infty,n]\cap \Z$ coincides with the involution $(i,i+1)$.
\end{lemma}

\begin{proof}
Note that $\sigma$ is a product of pairwise disjoint transpositions of the form $(m,m+1)$.
We fix $i$ and proceed by induction on $n$ satisfying that $n>i$.
For the base case, $n=i+1$, the element $z^{-i}\sigma z^i$ suffices.
Now assume that the element $\nu$ has been chosen so that $\nu\restriction_{(-\infty,n]\cap \Z}=(i,i+1)$.
Now $(n+1)\cdot \nu=m$ satisfies that $m\geq n+1$. If $m=n+1$, then we are done.
If $m>n+1$, we consider the element $f=\nu z^{-(m-1)}\sigma z^{(m-1)}$.
This element satisfies that $f\restriction_{(-\infty,n]\cap \Z}=(i,i+1)$ and that $(n+1)\cdot f=m-1$.
If $n+1=m-1$, then $f$ is the required element. Otherwise, we consider the element $g=fz^{-(m-2)}\sigma z^{(m-2)} $ which satisfies that
$g\restriction (-\infty,n]\cap \Z=(i,i+1)$ and that $(n+1)\cdot g=m-2$.
Proceeding in this manner, step by step, we obtain an element with the desired properties.
\end{proof}

\begin{prop}
The group $G$ is highly transitive and MIF.
\end{prop}

\begin{proof}
It is an immediate consequence of Lemma \ref{lem:HT} that for each $n\in \N$, there is a subgroup $H\leq G$ such that $P=\{-n,\ldots, n\}$ is $H$ invariant and $H\restriction P$ agrees with $\Sym(P)$. It follows that $G$ acts highly transitively on $\Z$.
To show that it is MIF, it suffices to show that $G\cap \Sym_f(\Z)=\{1\}$, using the Hull-Osin dichotomy stated in the introduction.
Since $\sigma$ is an involution, any element of the group $G$ (in the generating set $\{z,\sigma\}$) can be represented as a word of the form: $$z^{m_0} (z^{-m_1} \sigma z^{m_1})\ldots (z^{-m_k} \sigma z^{m_k})$$
Since $\sup\{m_0,\ldots, m_k\}<\infty$ and since $$\lim_{n\to \infty} |\inf(\Supp(\sigma_{n+1}\sigma_n^{-1}))-\sup(\Supp(\sigma_n))|\to \infty$$
we obtain that any nontrivial element of $G$ must have infinite support when viewed as a permutation in $\Sym(\Z)$.
\end{proof}

\bibliographystyle{alpha}
\bibliography{MIF.bib}

\end{document}